\documentclass[11pt,leqno]{amsart}

\usepackage{graphicx}
\usepackage{geometry}
\usepackage{amsmath}
\usepackage{amssymb}
\usepackage{amsthm}
\usepackage{amsfonts}
\usepackage{enumerate}
\usepackage[scriptsize,hang,raggedright]{subfigure}
\usepackage[cmyk]{xcolor}

\usepackage{hyperref}
\usepackage{pdfsync}
\usepackage{dsfont}
\usepackage{color}

\numberwithin{equation}{section} 

\newtheorem{theorem}{Theorem}[section]
\newtheorem{proposition}[theorem]{Proposition}
\newtheorem{lemma}[theorem]{Lemma}

\theoremstyle{remark}
\newtheorem{remark}[theorem]{Remark}

\newcommand{\R}{\mathbb{R}}
\newcommand{\disp}{\displaystyle}

\newcommand{\ba}{\begin{array}}
\newcommand{\ea}{\end{array}}

\newcommand{\bthm}{\begin{theorem}}
\newcommand{\ethm}{\end{theorem}}
\newcommand{\bprop}{\begin{proposition}}
\newcommand{\eprop}{\end{proposition}}
\newcommand{\blemma}{\begin{lemma}}
\newcommand{\elemma}{\end{lemma}}

\newcommand{\beqn}{\begin{equation}}
\newcommand{\eeqn}{\end{equation}}
\newcommand{\beqns}{\begin{equation*}}
\newcommand{\eeqns}{\end{equation*}}

\newcommand{\supp}{\operatorname{supp}}

\newcommand{\pr}{\prime}
\newcommand{\pt}{\partial}
\newcommand{\arrow}{\rightarrow}

\newcommand{\rn}{\mathbb{R}^N}

\newcommand{\Prob}{\mathcal{P}}

\renewcommand{\leq}{\leqslant}
\renewcommand{\geq}{\geqslant}

\newcommand{\te}{\textrm}

\newcommand{\veps}{\varepsilon}
\newcommand{\rhu}{\rightharpoonup}

\DeclareMathOperator{\divergence}{div}

\definecolor{mygreen}{rgb}{0.1,0.75,0.2}

\title[Existence of Ground States of Nonlocal-Interaction Energies]{Existence of Ground States of Nonlocal-Interaction Energies}

\author{Robert Simione}
\address{Department of Mathematical Sciences, Carnegie Mellon University, Pittsburgh, PA USA}
\email{rsimione@andrew.cmu.edu}
\author{Dejan Slep\v{c}ev}
\address{Department of Mathematical Sciences, Carnegie Mellon University, Pittsburgh, PA USA}
\email{slepcev@math.cmu.edu}
\author{Ihsan Topaloglu}
\address{Department of Mathematics and Statistics, 
McMaster University,
Hamilton, ON Canada}
\email{ihsan.topaloglu@math.mcmaster.ca}

\date{\today}                                        

\begin{document}

\begin{abstract}
We investigate  which nonlocal-interaction energies have a ground state (global minimizer). 
We consider this question over the space of probability measures and
 establish a sharp condition for the existence of ground states.  We show that this condition is closely related to the notion of stability (i.e. $H$-stability) of pairwise interaction potentials. Our approach uses the direct method of the calculus of variations. 
\end{abstract}

\maketitle

\section{Introduction}\label{sec:intro}

We investigate the existence of ground states (global minimizers) of nonlocal-interaction energies
	\beqn
		E(\mu) := \int_{\rn}\!\int_{\rn} w(x-y)\,d\mu(x)d\mu(y)
		\label{eqn:energy}
	\eeqn
considered over the space of probability measures $\Prob(\rn)$. 
Nonlocal-interaction energies arise in descriptions of systems of interacting particles, as well as their continuum limits. They are important to statistical mechanics \cite{FishRuel, Ruelle, Suto},  models of collective behavior of many-agent systems 
\cite{BeTo2011, M&K},  granular media \cite{Benedetto_etal,CaMcVi2006,Toscani2000}, self-assembly of nanoparticles \cite{HoPu2005,HoPu2006},  crystallization  \cite{AuFrSc, Radin, Theil}, and molecular dynamics simulations of matter \cite{Haile1992}.

Whether the energy dissipated by a system admits a global minimizer has important consequences on the behavior of the system. Continuum  systems governed by the energy which has a ground state typically exhibit  well defined dense clumps, while the systems with no global minimizers tend to disperse indefinitely. 

The interaction potential $w$ depends on the system considered. In most cases
it depends only on the distance between particles/agents. That is
the interaction potential $w$ is radially symmetric: $w(x)=W(|x|)$ for some $W:[0,\infty)\arrow (-\infty, \infty]$. Many potentials considered in the applications are repulsive at short distances ($W^\pr(r)<0$ for $r$ small) and attractive at large distances ($W^\pr(r)>0$ for $r$ large). 
Systems with finitely many particles governed by short-range-repulsive, long-range-attractive 
interaction potentials form well defined structures  (crystals are an example \cite{Theil}).  
The relevance of our result is to the behavior of these systems as the number of particles grows to infinity. Systems which have a global minimizer over the space of measures  form well defined states whose density grows as the number of particles increases, while the systems with no ground states typically  have bounded density and increase in size indefinitely.  

This mirrors the considerations in classical statistical mechanics when thermodynamic limit of particle systems is considered \cite{Ruelle}. Here we obtain mathematical results that highlight the connection. Namely, in Theorem \ref{thm:existence-vanish} (combined with Proposition \ref{ruelles}), we establish that the sharp condition for the existence of ground states of \eqref{eqn:energy} is closely related to the notion of stability ($H$-stability) of interacting potentials \cite{FishRuel, Ruelle}. More precisely we show that systems admitting a minimizer of  \eqref{eqn:energy} are (almost) precisely those for which the interaction potential is not $H$-stable, that is those for which the potential is catastrophic.

\medskip

In recent years significant interest in nonlocal-interaction energies arose from studies of
dynamical models.
For semi-convex  interaction potentials $w$ a number of systems governed by the energy $E$ can be interpreted as a gradient flow of the energy with respect to the Wasserstein metric and satisfy the
 nonlocal-interaction equation
	\beqn
		\frac{\pt\mu}{\pt t}=2 \divergence\left(\mu(\nabla w * \mu)\right).
		\label{eqn:cont-eqn}
	\eeqn
 Applications of this equation include models of collective behavior in biology \cite{BeTo2011, M&K},  granular media \cite{Benedetto_etal,CaMcVi2006,Toscani2000}, and self-assembly of nanoparticles \cite{HoPu2005,HoPu2006}.

While  purely attractive potentials lead to finite-time or infinite time blow up \cite{BertozziCarilloLaurent} the attractive-repulsive  potentials often generate finite-sized, confined aggregations \cite{FeHuKo11, KoSuUmBe2011, LeToBe2009}. 
The study of the nonlocal-interaction equation \eqref{eqn:cont-eqn} in terms of well-posedness, finite or infinite time blow-up, and long-time behavior has attracted the interest of many research groups in the recent years \cite{Balague_etal12,BaCaYa,BertozziCarilloLaurent,BertozziLaurent,BertozziLaurentLeger,BertozziLaurentRosado,CaDiFiLaSl2011,CaDiFiLaSl2012,FeRa10,FeHuKo11,KoHuPa13,KoSuUmBe2011,Laurent2007}.
The energy \eqref{eqn:energy} plays an important role in these studies as it governs the dynamics and as its (local) minima describe the long-time asymptotics of solutions. 

It has been observed that even for quite simple repulsive--attractive potentials the ground states are sensitive to the precise form of the potential and can exhibit a wide variety of patterns \cite{KoHuPa13,KoSuUmBe2011, vBUKB}.
In \cite{Balague_etal13}  Balagu\'e, Carrillo,  Laurent, and Raoul obtain conditions for the dimensionality of the support of local minimizers of \eqref{eqn:energy} in terms of the repulsive strength of the potential $w$ at the origin. Properties of steady states for a special class of potentials which blow up approximately like the Newtonian potential  at the origin have also been studied  \cite{BertozziLaurentLeger,  CDM14,  FeHu13, FeHuKo11}.
Particularly relevant to our study are the results obtained by Choksi, Fetecau and one of the authors \cite{ChFeTo14} on the existence of minimizers of interaction energies in a certain form. There the authors consider potentials of the power-law form, $w(x):=|x|^a/a-|x|^r/r$, for $-N<r<a$, and prove the existence of minimizers in the class of probability measures when the power of repulsion $r$ is positive. When the interaction potential has a singularity at the origin, i.e., for $r<0$, on the other hand, they establish the existence of minimizers of the interaction energy in a restrictive class of uniformly bounded, radially symmetric $L^1$-densities satisfying a given mass constraint.  Carrillo, Chipot and Huang \cite{CaChHu} also consider the minimization of nonlocal-interaction energies defined via power-law potentials and prove the existence of a global minimizer by using a discrete to continuum approach. The ground states and their relevance to statistical mechanics were also considered in periodic setting (and on bounded sets) by S\"uto \cite{Suto}.

\subsection{Outline.}
In Theorems \ref{thm:existence-blow-up} and \ref{thm:existence-vanish} we establish criteria for the existence of minimizers of a very broad class of potentials. We employ the direct method of the calculus of variations. In Lemma \ref{LSCE} we establish the weak lower-semicontinuity of the energy with respect to weak convergence of measures. When the potential $W$ grows unbounded at infinity  (case treated in Theorem \ref{thm:existence-blow-up}) this provides enough confinement for a minimizing sequence to ensure the existence of minimizers. If $W$ asymptotes to a finite value  (case treated in Theorem \ref{thm:existence-vanish}) then there is a delicate interplay between repulsion at some lengths (in most applications short lengths) and attraction at other length scales (typically long) which establishes whether 
the repulsion wins and a minimizing sequence spreads out indefinitely and ``vanishes'' or the minimizing sequence is compact and has a limit. We establish a simple, sharp condition, (\textbf{HE}) on the energy that characterizes whether a ground state exists.
To establish compactness of a minimizing sequence we use Lions' concentration compactness lemma.


While the conditions (\textbf{H1}) and (\textbf{H2}) are easy-to-check conditions on the potential $W$ itself, the condition  (\textbf{HE}) is a condition on the energy and it is not always easy to verify. Due to the above connection with statistical mechanics the conditions on $H$-stability (or the lack thereof) can be used to verify if (\textbf{HE}) is satisfied for a particular potential. We list such conditions in Section \ref{sec:cond-on-w}. However only few general conditions are available. It is an important open problem to establish a more complete characterization of potentials $W$ which satisfy (\textbf{HE}). 

We finally remark that as this manuscript was being completed we learned that  Ca\~nizo, Carrillo, and  Patacchini \cite{CCP} independently and concurrently obtained very similar conditions for the existence of minimizers, which they also show to be compactly supported. 
The proofs however are quite different. 

\section{Hypotheses and Preliminaries}\label{sec:hypo}
	
The interaction potentials we consider are radially symmetric, that is, $w(x)=W(|x|)$ for some function $W:[0,\infty)\arrow \mathbb{R}\cup\{\infty\}$, and they satisfy the following basic properties:
\medskip

\begin{itemize}	
\addtolength{\itemsep}{3pt}
\item[({\bf H1})] $W$ is lower-semicontinuous. 
\item[({\bf H2})] The function $w(x)$ is locally integrable on $\R^N$.
\end{itemize}
	
\medskip

Beyond the basic assumptions above, the behavior of the tail of $W$ will play an important role.
We consider potentials which have a limit at infinity. If the limit is finite we can add a constant to the potential, which does not affect the existence of minimizers, and assume that the limit is zero. 
 If the limit is infinite the proof of existence of minimizers is simpler, while if the limit is finite an additional condition is needed. Thus we split the condition on behavior at infinity 
into two conditions:
\medskip

	\begin{itemize}
	\addtolength{\itemsep}{3pt}
		\item[({\bf H3a})] $W(r) \arrow \infty$ as $r\arrow\infty$.
	\end{itemize}
  \begin{itemize}
		\item[({\bf H3b})] $W(r)\arrow 0$ as $r\arrow\infty$.
	\end{itemize}

\medskip

\begin{remark}\label{rmk:boundedbelow}
By the assumptions ({\bf H1}) and ({\bf H3a}) or ({\bf H3b}) the interaction potential $W$ is bounded from below. Hence
\begin{equation} \label{CW}
 C_W :=  \inf_{r \in (0, \infty)} W(r) > -\infty . 
 \end{equation} 
If ({\bf H3a}) holds,  by adding $-C_W$ to $W$ from now on we assume that $W(r) \geq 0$ for all $r \in (0, \infty)$
\end{remark}

\begin{figure}[ht!]
     \begin{center}

        \subfigure[Interaction potentials satisfying ({\bf H1}),  ({\bf H2}), and ({\bf H3a})]{
            \label{fig:first}         \includegraphics[width=0.3\linewidth]{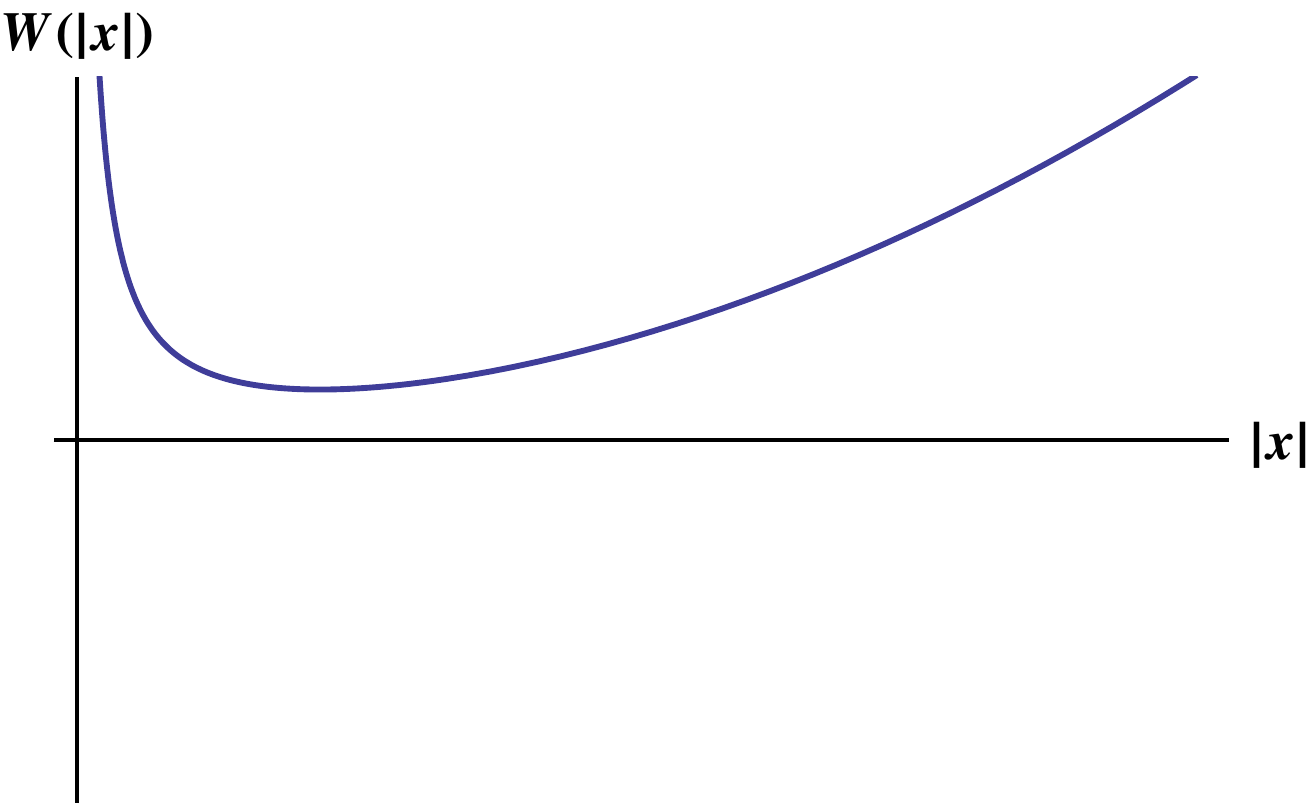}\qquad\qquad\qquad\includegraphics[width=0.3\linewidth]{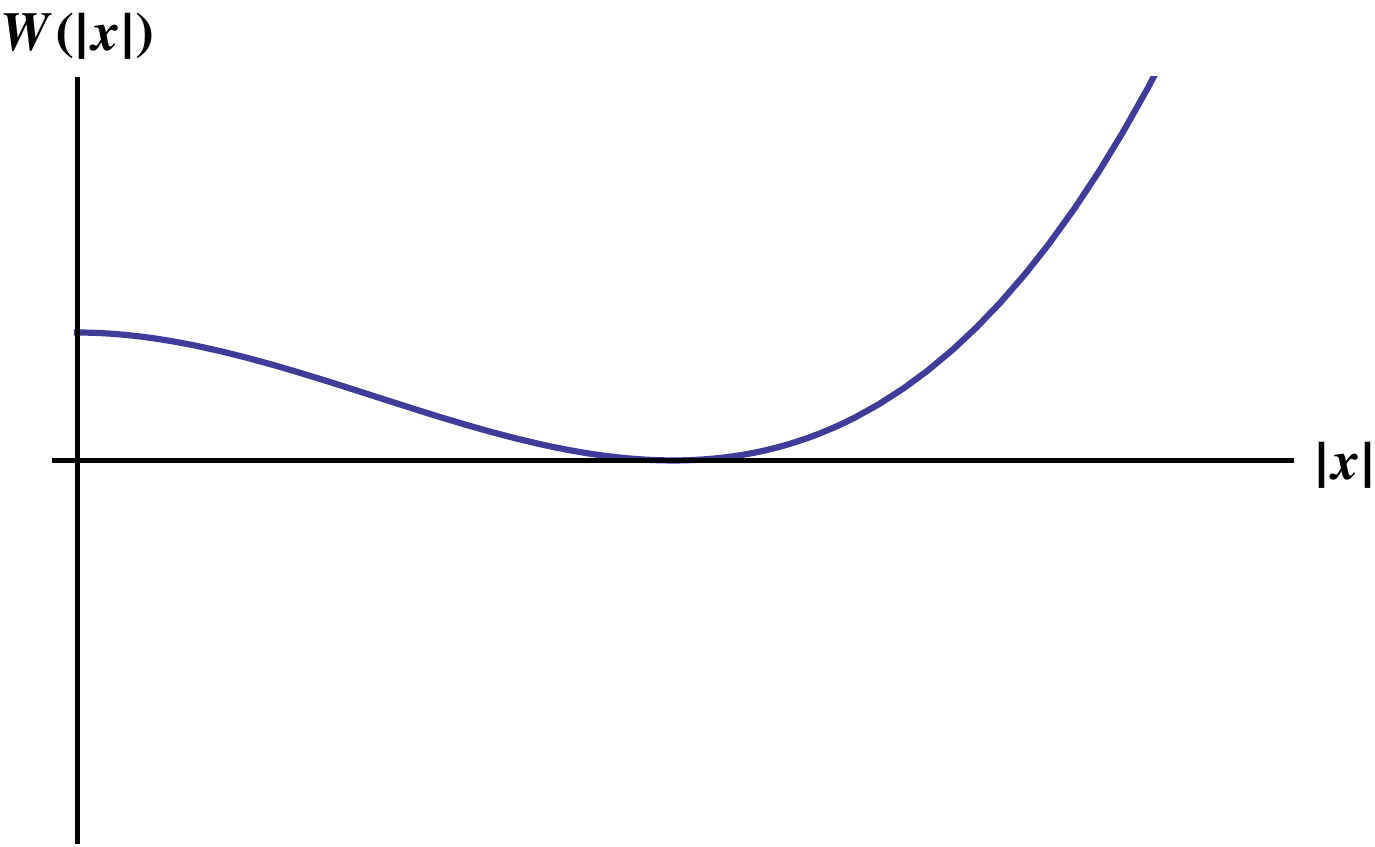}
        }
        \\ 
        \subfigure[Interaction potentials satisfying ({\bf H1}), ({\bf H2}), and ({\bf H3b})]{
            \label{fig:second}       \includegraphics[width=0.3\linewidth]{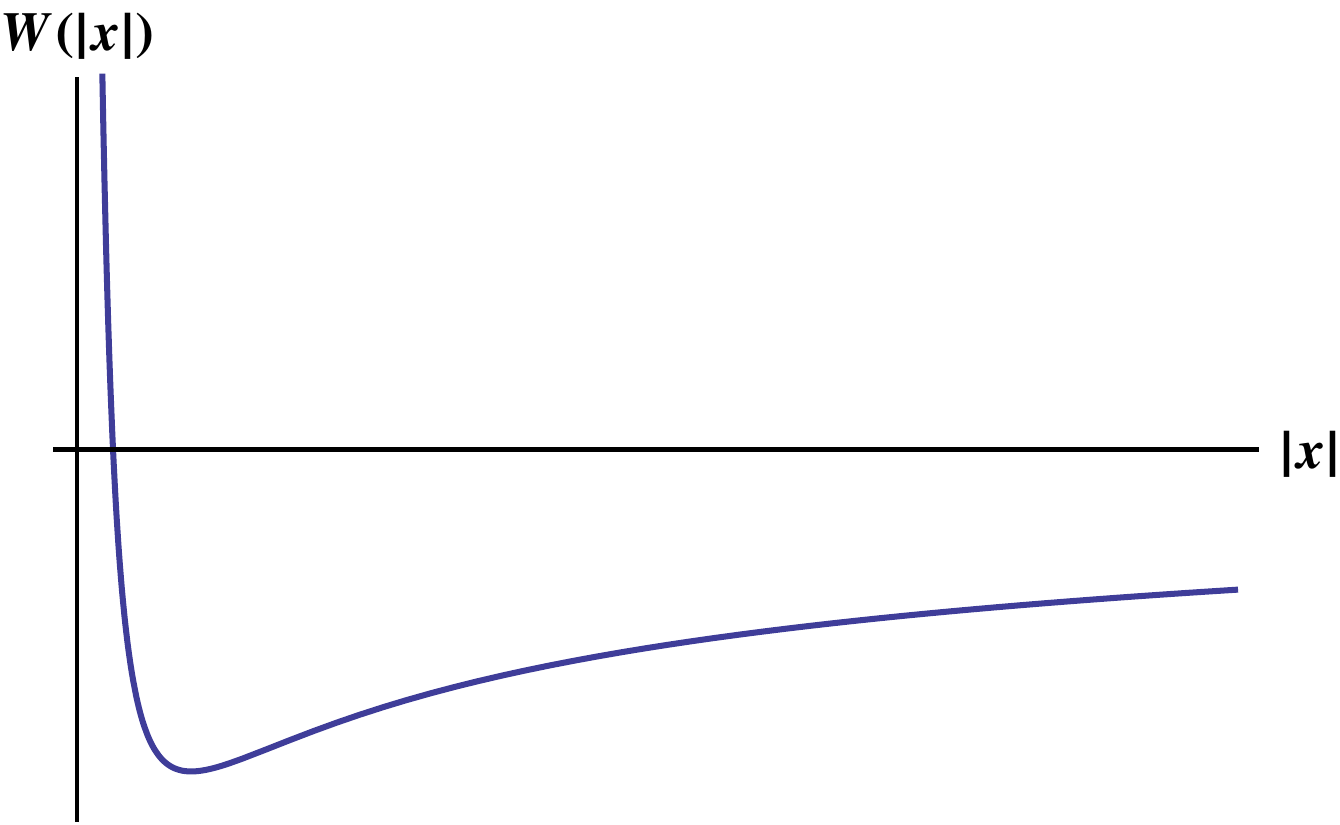}\qquad\qquad\qquad\includegraphics[width=0.3\linewidth]{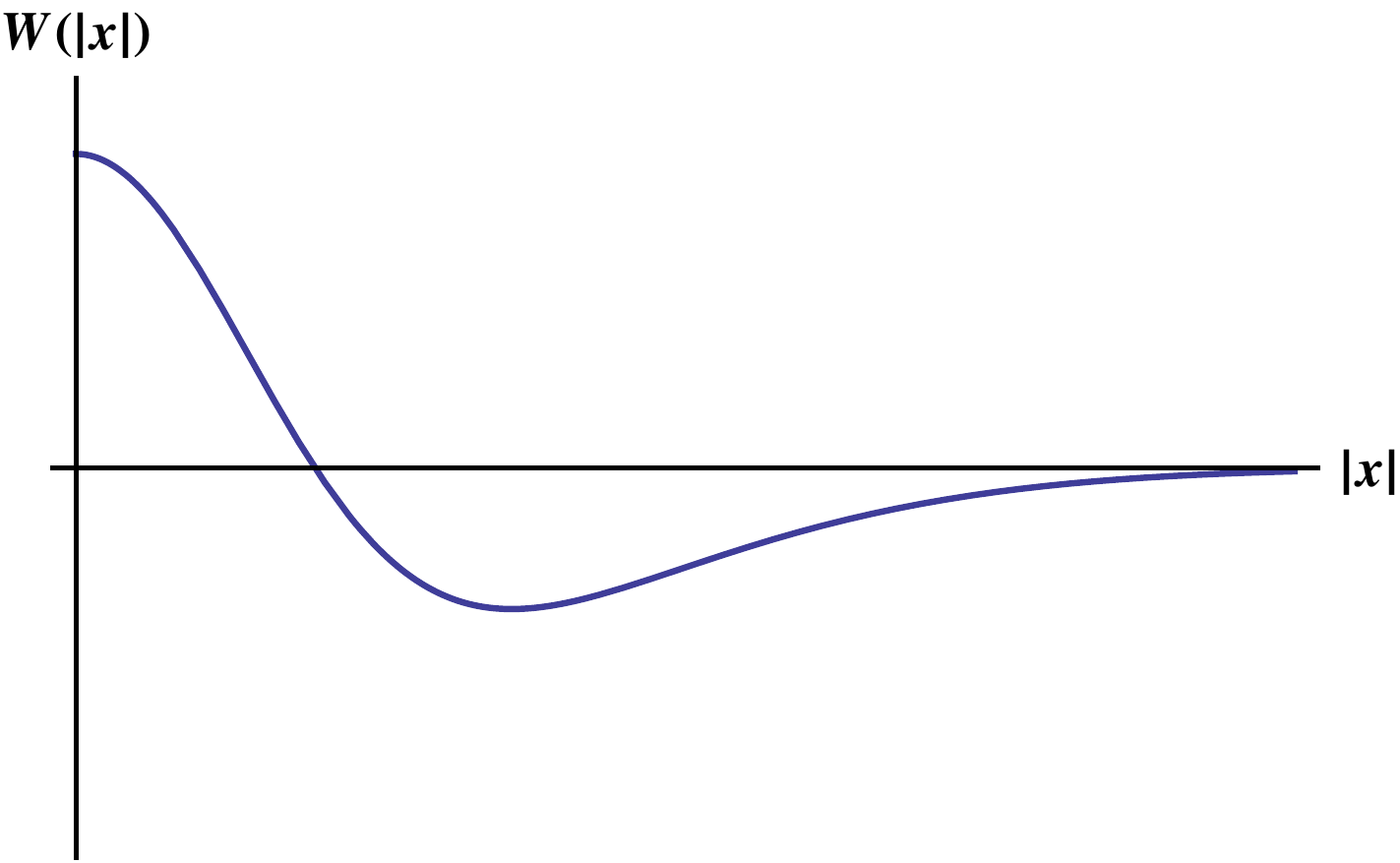}
        }
    \end{center}
    \caption{
        Generic examples of $W(|x|)$.
     }
   \label{fig:potential-examples}
\end{figure}

As noted in the introduction the assumptions ({\bf H1}), ({\bf H2}) with ({\bf H3a}) or ({\bf H3b}) allow us to handle a quite general class of interaction potentials $w$.  Figure \ref{fig:potential-examples} illustrates a set of simple examples of smooth potential profiles $W$ that satisfy these assumptions.

\medskip

In order to establish the existence of ground states of $E$, for interaction potentials $w$ satisfying ({\bf H1}), ({\bf H2}) and ({\bf H3b}), the  following assumption on the interaction energy $E$ is needed:
\medskip

({\bf HE})  There exists a measure $\bar{\mu}\in\Prob(\rn)$ such that $E(\bar{\mu}) \leq 0$.
\medskip

We establish that the conditions ({\bf H1}), ({\bf H2}) and ({\bf H3a}) or ({\bf H3b})  imply the lower-semicontinuity of the energy with respect to weak convergence of measures. We recall that a sequence of probability measures $\mu_n$ converges 
weakly to measure $\mu$, and we write $\mu_n \rhu \mu$, if for every bounded continuous function $\phi \in C_b(\R^N,\R)$
\[ \int \phi d \mu_n \to \int \phi d \mu \quad \te{ as } n \to \infty. \]

\blemma[Lower-semicontinuity of the energy] \label{LSCE}
Assume $W : [0, \infty) \to (-\infty, \infty]$ is a lower-semicontinuous function bounded from below. 
Then the energy $E: \mathcal P(\R^n) \to (-\infty, \infty]$ defined in \eqref{eqn:energy} is weakly lower-semicontinuous with respect to weak convergence of measures.
\elemma
\begin{proof}
Let $\mu_n$ be a sequence of probability measures such that $\mu_n \rhu \mu$ as $n \to \infty$.
Then $\mu_n \times \mu_n \rhu \mu \times \mu$ in the set of probability measures on $\R^N \times \R^N$. If $w$ is  continuous and bounded 
\[ \int_{\rn}\!\int_{\rn} w(x-y)\,d\mu_n(x)d\mu_n(y) \longrightarrow \int_{\rn}\!\int_{\rn} w(x-y)\,d\mu(x)d\mu(y) \quad \te{ as } n \to \infty. \]
So, in fact, the energy is continuous with respect to weak convergence. 
On the other hand, if $w$ is lower-semicontinuous and $w$ is bounded from below then the weak lower-semicontinuity of the energy follows from the Portmanteau Theorem \cite[Theorem 1.3.4]{van1996weak}.
\end{proof}
We remark that the assumption on boundedness from below is needed since if, for example, $W(r) = -r$ then for $\mu_n = (1-\frac{1}{n}) \delta_0 + \frac{1}{n} \delta_n$ the energy is $E(\mu_n) = -1$ for all $n\in\mathbb{N}$, while $\mu_n \rhu \delta_0$ which has energy $E(\delta_0)=0$.
\medskip

Finally, we state Lions' concentration compactness lemma for probability measures
\cite{Lions84}, \cite[Section 4.3]{Struwe}.  We use this lemma to verify that an energy-minimizing sequence is precompact in the sense of weak convergence of measures.

\blemma[Concentration-compactness lemma for measures]
		Let $\{\mu_n\}_{n\in\mathbb{N}}$ be a sequence of probability measures on $\rn$. Then there exists a subsequence $\{\mu_{n_k}\}_{k\in\mathbb{N}}$ satisfying one of the three following possibilities:
		\vspace{0.2cm}
	\begin{itemize}
		\item[(i)] \emph{(tightness up to translation)} There exists a sequence $\{y_k\}_{k\in\mathbb{N}}\subset\rn$ such that for all $\veps>0$ there exists $R>0$ with the property that
			\[
				\int_{B_R(y_k)}\,d\mu_{n_k}(x) \geq 1-\veps \qquad {\hbox{\rm  for all $k$.}}
			\]
		 \item[(ii)] \emph{(vanishing)} $\disp \lim_{k\arrow\infty} \sup_{y\in\rn} \int_{B_R(y)}\,d\mu_{n_k}(x)=0$, for all $R>0$; \vspace{0.2cm}
		 \item[(iii)] \emph{(dichotomy)} There exists $\alpha\in(0,1)$ such that for all $\veps>0$, there exist a number $R>0$ and a sequence $\{x_k\}_{k\in\mathbb{N}}\subset\rn$ with the following property:\\
		 
		 	\vspace{-0.2cm}
		 Given any $R^\pr>R$ there are nonnegative measures $\mu_k^1$ and $\mu_k^2$ such that
		 	\vspace{0.2cm}
			\begin{itemize}
		 		\item[] $0\leq \mu_k^1 + \mu_k^2 \leq \mu_{n_k}$\,,
				\vspace{0.2cm}
		 		\item[] $\supp(\mu_k^1)\subset B_R(x_k)$,\qquad $\supp(\mu_k^2)\subset \rn\setminus B_{R^\pr}(x_k)$\,,
				\vspace{0.2cm}
		 		\item[] $\disp \limsup_{k\arrow\infty} \left(\left|\alpha-\int_{\rn}d\mu_k^1(x)\right|+\left|(1-\alpha)-\int_{\rn}d\mu_k^2(x)\right|\right)\leq \veps$.
		 	\end{itemize}
	\end{itemize}
	\label{cmptlem}
\elemma

%

\section{Existence of Minimizers}\label{sec:existence}

In this section we prove the existence of a global minimizer of $E$. We use the direct method of the calculus of variations and utilize Lemma \ref{cmptlem} to eliminate the ``vanishing'' and ``dichotomy'' of an energy-minimizing sequence. The techniques in our proofs, though, depends on the behavior of the interaction potential at infinity. Thus we prove two existence theorems: one for potentials satisfying ({\bf H3a}) and another one for those satisfying ({\bf H3b}).

\bthm\label{thm:existence-blow-up}
	Suppose $W$ satisfies the assumptions ({\bf H1}), ({\bf H2}) and ({\bf H3a}). Then the energy \eqref{eqn:energy} admits a global minimizer in $\Prob(\rn)$.
\ethm

\begin{proof}
 Let $\{\mu_n\}_{n\in\mathbb{N}}$ be a minimizing sequence, that is, $\lim_{n\arrow\infty} E(\mu_n)=\inf_{\mu\in\Prob(\rn)} E(\mu)$.

\medskip
Suppose $\{\mu_k\}_{k\in\mathbb{N}}$ has a subsequence which ``vanishes''. Since that subsequence is also a minimizing sequence we can assume that $\{\mu_k\}_{k\in\mathbb{N}}$ vanishes.
  Then  for any $\veps>0$ and for any $R>0$ there exists $K\in\mathbb{N}$ such that for all $k>K$ and for all $x\in\rn$
 	\[
 		\mu_{k}(\rn\setminus B_R(x))\geq 1-\veps.
 	\]
This implies that for $k>K$,
	\[
		\iint_{|x-y|\geq R} d\mu_{k}(x)d\mu_{k}(y) = \int_{\rn}\!\left(\int_{\rn\setminus B_R(x)} d\mu_{k}(y)\right)\,d\mu_{k}(x) \geq 1-\veps.
	\]
Given $M \in \R$, by condition (\textbf{H3a}) there exists $R>0$ such that for all $r \geq R$, $W(r) \geq M$.
Consider $\veps \in (0, \frac12)$ and $K$ corresponding to $\veps$ and $R$.
Since $W \geq 0$ by Remark \ref{rmk:boundedbelow},
	\beqns
			\begin{aligned}
				E(\mu_{k}) &= \iint_{|x-y| < R} W(|x-y|)\,d\mu_{k}(x)d\mu_{k}(y) + \iint_{|x-y|\geq R} W(|x-y|)\,d\mu_{k}(x)d\mu_{k}(y) \\
										 &\geq \iint_{|x-y|\geq R} W(|x-y|)\,d\mu_{k}(x)d\mu_{k}(y) \\
										 &\geq (1-\veps)M
			\end{aligned}
	\eeqns
for all $k>K$. Letting $M \to \infty$ implies $E(\mu_{k})\arrow\infty$. This contradicts the fact that $\mu_{k}$ is a subsequence of a minimizing sequence of $E$. Thus, ``vanishing'' does not occur.

Next we show that ``dichotomy'' is also not an option for a minimizing sequence. Suppose,  that ``dichotomy'' occurs.  As before we can assume that the subsequence along which dichotomy occurs is the whole sequence. Let $\veps>0$ be fixed, and let $R$, the sequence $\{x_k\}_{k\in\mathbb{N}}$ and measures
	\[
		\mu_k^1+\mu_k^2 \leq \mu_{k}.
	\]
be as defined in Lemma \ref{cmptlem}(ii).
For any $R'>R$	, using Remark \ref{rmk:boundedbelow}, we obtain
	\beqns
		\begin{aligned}
		\liminf_{k\arrow\infty} E(\mu_{n_k}) & \geq  \liminf_{k\arrow\infty} \int_{B_R(x_{n_k})}\!\int_{B^c_{R^\pr}(x_{n_k})} W(|x-y|)\,d\mu_k^2(x)d\mu_k^1(y)\\
								 &\geq  \inf_{r \geq R'-R}  W(r)\, (\alpha-\veps)(1-\alpha-\veps) ,
		\end{aligned}
	\eeqns
where $B^c_{R^\pr}(x_{n_k})$ simply denotes $\rn\setminus B_{R^\pr}(x_{n_k})$.

By ({\bf H3a}), letting $R^\pr\arrow\infty$ yields that
	\[
		\liminf_{k\arrow\infty} E(\mu_{n_k}) \geq \infty,
	\]
which contradicts the fact that $\mu_k$ is an energy minimizing sequence.

Therefore ``tightness up to translation'' is the only possibility. Hence there exists $y_k\in\rn$ such that for all $\veps>0$ there exists $R>0$ with the property that
			\[
				\int_{B(y_k,R)}\,d\mu_{n_k}(x) \geq 1-\veps \qquad {\hbox{\rm  for all $k$.}}
			\]
Let
	\[
		\tilde{\mu}_{n_k}:=\mu_{n_k}(\cdot + y_k).
	\]
Then the sequence of probability measures $\{\tilde{\mu}_{n_k}\}_{k\in\mathbb{N}}$ is tight.
Since the interaction energy is translation invariant we have that
	\[
		E(\tilde{\mu}_{n_k})=E(\mu_{n_k}).
	\]
Hence, $\{\tilde\mu_{n_k}\}_{k\in\mathbb{N}}$ is also an energy-minimizing sequence. 
By the Prokhorov's theorem (cf. \cite[Theorem 4.1]{Billingsley}) there exists a further subsequence of $\{\tilde{\mu}_{n_k}\}_{k\in\mathbb{N}}$ which we still index by $k$, and a measure $\mu_0\in\Prob(\rn)$ such that
	\[
		\tilde{\mu}_{n_k} {\rightharpoonup}\mu_0
	\]
in $\mathcal{P}(\rn)$ as $k\arrow\infty$.

\medskip

Since the energy in lower-semicontinuous with respect to weak convergence of measures, by Lemma
\ref{LSCE}, the measure $\mu_0$ is  a minimizer of $E$.
\end{proof}

The second existence theorem involves interaction potentials which vanish at infinity. 

\bthm\label{thm:existence-vanish}
Suppose $W$ satisfies the assumptions ({\bf H1}), ({\bf H2}) and ({\bf H3b}).  Then the energy $E$,  given by \eqref{eqn:energy}, has a global minimizer in $\Prob(\rn)$ if and only if it satisfies the condition ({\bf HE}).
\ethm

\begin{proof}
Let us assume that $E$ satisfies  condition ({\bf HE}).
 As before, our proof relies on the direct method of the calculus variations for which we need to establish precompactness of a minimizing sequence.

Let $\{\mu_n\}_{n\in\mathbb{N}}$ be a minimizing sequence and let
	\[
		I:=\inf_{\mu\in\Prob(\rn)} E(\mu).
	\]
Condition  ({\bf HE}) implies that  $I \leq 0$.
If $I = 0$ then by assumption ({\bf HE}) there exists $\bar \mu$ with $E(\bar\mu)=0$, which is
the desired minimizer.  Thus, we focus on case that $I<0$.  Hence there exists $\bar \mu$ for which 
 $E(\bar{\mu})<0$. Also note that by Remark \ref{rmk:boundedbelow}, $I>-\infty$.

\medskip

Suppose the subsequence $\{\mu_{n_k}\}_{k\in\mathbb{N}}$ of the minimizing sequence $\{\mu_n\}_{n\in\mathbb{N}}$ ``vanishes''.  Since that subsequence is also a minimizing sequence we can assume that $\{\mu_k\}_{k\in\mathbb{N}}$ vanishes. That is, for any $R>0$ 
	\beqn
		\lim_{k\arrow\infty} \sup_{x \in \R^N} \int_{B_R(x)} d\mu_{k}(y) = 0.
		\label{eqn:vanishing-second-exist}
	\eeqn
Let \[\overline W(R) = \inf_{r \geq R} W(r). \]
 Since $W(r) \to 0$ as $r \to \infty$, $\overline W(r) \to 0$ as $r \to \infty$ and $\overline W(r) \leq 0$ for all $r \geq 0$. Then we have that  
	\beqns
		\begin{aligned}
			E(\mu_{k}) &= \iint_{|x-y| > R} W(|x-y|)\,d\mu_{k}(x)d\mu_{k}(y) + \iint_{|x-y|\leq R} W(|x-y|)\,d\mu_{k}(x)d\mu_{k}(y) \\
					 & \geq \overline W(R) + C_W \iint_{|x-y|\leq R} d\mu_{k}(x)d\mu_{k}(y) \\
					&= \overline W(R) + C_W \int_{\rn}\!\left(\int_{B_R(x)}d\mu_{k}(y)\right)\,d\mu_{k}(x).
		\end{aligned}
	\eeqns
Vanishing of the measures, \eqref{eqn:vanishing-second-exist}, implies that
	$		\liminf_{k\arrow\infty} E(\mu_{k}) \geq \overline W(R)
	$
for all $R>0$. Taking the limit as $R \to \infty$ gives 
	\[
		\liminf_{k\arrow\infty} E(\mu_{k}) \geq 0.
	\]
This contradicts the fact that the infimum of the energy, namely $I$, is negative. Therefore ``vanishing'' in Lemma \ref{cmptlem} does not occur.
\medskip

Suppose the dichotomy occurs.
Let $\alpha\in(0,1)$  and $R>0$ be as in Lemma \ref{cmptlem}  and $C_W$ be the constant defined in \eqref{CW}.
Let $\veps >0$ be such that 
\begin{equation} \label{epsbound}
\veps < \frac{|I|}{64 |C_W|} \min \left\{\frac{1}{\alpha}-1, \frac{1}{1-\alpha}-1 \right\}
\end{equation}
and let $R'$ be such that 
\begin{equation} \label{Rbound}
|\overline W(R'-R)| = |\inf_{r \geq R'-R} W(r)| < \frac{|I|}{32} \min \left\{\frac{1}{\alpha}-1, \frac{1}{1-\alpha}-1 \right\}.
\end{equation} 
 As in the proof of Theorem \ref{thm:existence-blow-up}, we can assume that  dichotomy occurs along the whole sequence. 
Let $ \mu_k^1$ and  $\mu_k^2$ be measures described in Lemma \ref{cmptlem}. Let $\nu_k =  \mu_k -  (\mu_k^1 + \mu_k^2)$.  Note that $\nu_k$ is a nonnegative measure with $|\nu_k|< \veps$, where 
$|\nu_k| = \nu_k(\R^N)$.

Let $B[\cdot,\cdot]$ denote the symmetric bilinear form
	\[
		B[\mu,\nu]:=2\int_{\rn}\!\int_{\rn} W(|x-y|)\,d\mu(x)d\nu(y).
	\]
By the definition of energy	
\begin{align}
\begin{split}
E(\mu_k) & = E(\mu_k^1) + E(\mu_k^2)+ B(\mu_k^1,\mu_k^2) + B(\mu_k^1+\mu_k^2, \nu_k) + E(\nu_k) \\
& \geq  	E(\mu_k^1) + E(\mu_k^2) - |\overline W(R'-R)| - 2|C_W| \veps
\end{split} \label{bee}
\end{align}
where we used that the supports of $\mu_k^1$ and $\mu_k^2$ are at least $R'-R$ apart. 
We can also assume, without loss of generality, that $E(\mu_k) < \frac12 I$ for all $k$.
Let $\alpha_k  = | \mu_k^1|$, $\beta_k=|\mu_k^2|$. 

Let us first consider the case that $\frac{1}{\alpha_k} E(\mu^1_k) \leq \frac{1}{\beta_k}  E(\mu^2_k)$. 
Note that the energy has the following scaling property:
	\[
		E(c\sigma)=c^2 E(\sigma)
	\]
for any constant $c>0$ and measure $\sigma$. Our goal is to show that for some $\lambda >0$, for all large enough $k$, $\,E(\frac{1}{\alpha_k} \mu_k^1 ) < E(\mu_k ) -  \lambda |I|$  which contradicts the fact that $\mu_k$ is a minimizing sequence.

Let us consider first the subcase that $E(\mu_k^2) \geq 0$ along a subsequence. By relabeling we can assume that the subsequence is the whole sequence.  From \eqref{epsbound}, \eqref{Rbound}, and  \eqref{bee} it follows  that
	\beqn\label{mu1bound}
		 \frac{1}{\alpha_k} E(\mu_k^1) < \frac{I}{4}
		\eeqn
for all $k$. 
Using the estimates again, we obtain
	\begin{align*}
 E(\mu_k) - E\left(\frac{1}{\alpha_k} \mu_k^1 \right) & \overset{\eqref{bee}}{\geq} \left(1 - \frac{1}{\alpha_k^2} \right) E\left(\mu_k^1 \right) - |\overline W(R'-R)| - 2|C_W| \veps \\
  &\overset{\eqref{mu1bound}}{ \geq} \left(\frac{1}{\alpha_k} -1 \right) \frac{|I|}{4} - |\overline W(R'-R)| - 2 |C_W| \veps \\
 &\!\!\!\!\!\overset{\eqref{epsbound},\eqref{Rbound}}{ \geq}  \left(\frac{1}{\alpha} -1 \right) \frac{|I|}{16}.
\end{align*}     
Thus $\mu_k$ is not a minimizing sequence. Contradiction.
 
 Let us now consider the subcase $E(\mu_k^2) \leq 0$ for all $k$. Using \eqref{bee} and $\frac{\beta_k}{\alpha_k} E(\mu^1_k) \leq  E(\mu^2_k)$ we obtain
 \[ \frac{I}{2} \geq E(\mu_k) \geq \left( 1 + \frac{\beta_k}{\alpha_k} \right) E(\mu_k^1) - |\overline W(R'-R)| -2|C_W| \veps. \]
 From \eqref{epsbound} and \eqref{Rbound} follows that for all $k$
 \[ \frac{1}{\alpha_k} E(\mu_k^1)  \leq \frac{I}{8}. \]
Combining with above inequalities gives
 \begin{align*}
 E(\mu_k) - E\left(\frac{1}{\alpha_k} \mu_k^1 \right) & \geq \left( 1 + \frac{\beta_k}{\alpha_k} - \frac{1}{\alpha_k^2}  \right) E(\mu_k^1) - | \overline W(R'-R)| - 2 |C_W| \veps \\
 & \geq  \left( \frac{1}{\alpha_k} - \alpha_k - \beta_k \right) \frac{|I|}{8} - \left( \frac{1}{\alpha} - 1 \right)  \left(\frac{|I|}{32} + \frac{|I|}{32} \right)\\
 & \geq \frac{|I|}{32} \left( \frac{1}{\alpha} -1 \right)
 \end{align*}
for $k$ large enough. This contradicts the assumption that $\mu_k$ is a minimizing sequence.

\medskip
The case $\frac{1}{\alpha_k} E(\mu^1_k) > \frac{1}{\beta_k}  E(\mu^2_k)$ is analogous.
In conclusion the dichotomy does not occur. Therefore ``tightness up to translation'' is the only possibility. As in the proof of Theorem \ref{thm:existence-blow-up}, 
we can translate measures $\mu_{n_k}$ to obtain a tight, energy-minimizing sequence $\tilde \mu_{n_k}$.

By Prokhorov's theorem, there exists a further subsequence of $\{\tilde\mu_{n_k}\}_{k\in\mathbb{N}}$, still indexed by $k$, such that
	\[
		\mu_{n_k} \rhu \mu_0 \quad \te{ as } k \to \infty
	\]
for some measure $\mu_0\in\Prob(\rn)$ in $\mathcal{P}(\rn)$ as $k\arrow\infty$.
Therefore, by lower-semicontinuity of the energy, $\mu_0$ is a minimizer of $E$ in the class $\Prob(\rn)$.
\medskip

We now show the necessity of  condition ({\bf HE}).
Assume that $E(\mu)> 0$ for all $\mu\in\Prob(\rn)$. To show that the energy $E$ does not have a minimizer  consider a sequence of measures which ``vanishes'' in the sense of Lemma \ref{cmptlem}(ii). Let
	\[
		\rho(x)=\frac{1}{\omega_N}\chi_{B_1(0)}(x),
	\]
where $\omega_N$ denotes the volume of the unit ball in $\rn$ and $\chi_{B_R(0)}$ denotes the characteristic function of $B_R(0)$, the ball of radius $R$ centered at the origin. Consider the sequence
	\[
		\rho_n(x)=\frac{1}{n^N}\rho\left(\frac{x}{n}\right)
	\]
for $n\geq 1$. Note that  $\rho_n$ are in $\Prob(\rn)$. We estimate
\beqns
	 	\begin{aligned}
		0<E(\rho_n)  &=	  \frac{1}{\omega_N^2 n^{2N}} \int_{B_n(0)}\!\int_{B_n(0)} W(|x-y|)\,dxdy\\
								 &\leq \frac{1}{\omega_N^2 n^{2N}} \int_{B_n(0)}\!\left(\int_{B_n(y)} |W(|x|)|\,dx\right)\,dy \\
								 &\leq \frac{1}{\omega_N n^{N}} \left( \int_{B_R(0)}|W(|x|)|\,dx+\int_{B_{2n}(0)\setminus B_R(0)}|W(|x|)|\,dx \right) \\
								 &\leq \frac{C(R)}{\omega_N n^{N}} +\frac{2^N}{\omega_N} \sup_{r \geq R} |W(r)|.  
		\end{aligned}
	\eeqns
Since $\sup_{r \geq R} |W(r)| \to 0$ as $R \to \infty$, for any $\veps>0$ we can choose $R$ so that
$\frac{2^N}{\omega_N} \sup_{r \geq R} |W(r)| < \frac{\veps}{2}$. We can then choose $n$ large enough for 
$\frac{C(R)}{\omega_N n^{N}} < \frac{\veps}{2}$ to hold. 
Therefore $\lim_{n\arrow\infty} E(\rho_n)=0$, that is, $\inf_{\mu\in\Prob(\rn)}E(\mu)=0$. However, since $E(\cdot)$ is positive for any measure in $\Prob(\rn)$ the energy does not have a minimizer.
\end{proof}

\section{Stability and Condition \textbf{(HE)}}\label{sec:cond-on-w}

The interaction energies of the form \eqref{eqn:energy} have been an important object of study in statistical mechanics. For a system of interacting particles to have a macroscopic thermodynamic behavior it is needed that it does not accumulate mass on bounded regions as the number of particles goes to infinity. Ruelle called such potentials stable (a.k.a. $H$-stable).
More precisely, a potential $W:[0,\infty) \to (-\infty, \infty]$ is defined to be \emph{stable} 
if there exists  $B \in \R$ such that for all $n$ and for all sets of $n$ distinct points $\{x_1, \dots, x_n\}$ in $\R^N$
\begin{equation} \label{Ruc}
\frac{1}{n^2} \sum_{1 \leq i < j \leq n} w(x_i-x_j) \geq - \frac{1}{n} B.
\end{equation}

We show that for a large class of pairwise interaction potentials the stability is equivalent with nonnegativity of energies. Our result is a continuum analogue of a part of \cite[Lemma 3.2.3]{Ruelle}.

\bprop[Stability conditions] \label{ruelles}
Let $W:[0, \infty) \to \R$ be an upper-semicontinuous function such that $W$ is bounded from above or
there exists $\overline R$ such that $W$ is nondecreasing on $[\overline R, \infty)$. Then the conditions 
\begin{itemize}
\item[(\textbf{S1})] $w$ is a stable potential as defined by \eqref{Ruc},
\item[(\textbf{S2})] for any probability measure $\mu \in \mathcal P(\R^N)$, $E(\mu) \geq 0$
\end{itemize}
are equivalent.
\eprop

Note that all potentials considered in the proposition are finite at $0$. We expect that the condition can be extended to a class of potentials which converge to infinity at zero. Doing so is an open problem.
We also note that the condition (\textbf{S2}) is not exactly the complement of (\textbf{HE}), as the nonnegative potentials whose minimum is zero satisfy both conditions. Such potentials indeed exist:
  for example consider any smooth nonnegative $W$  such that $W(0)=0$.
Then the associated energy is nonnegative and $E(\delta_0)=0$ so any singleton is an energy minimizer.
Note that $E$ satisfies both (\textbf{HE}) and stability. To further remark on connections with statistical mechanics we note that such potentials $W$ are not \emph{super-stable}, but are \emph{tempered}
if $W$ decays at infinity (both notions are defined in  \cite[Chapter 3]{Ruelle}).

\begin{proof}
To show that  (\textbf{S2}) implies (\textbf{S1}) consider $\mu = \frac{1}{n} \sum_{i=1}^n \delta_{x_i}$. Then from $E(\mu) \geq 0$ it follows that $\frac{1}{n^2} \sum_{1 \leq i < j \leq n} w(x_i-x_j) \geq - \frac{1}{2 n} W(0)$ so (\textbf{S1}) holds with  $B=\frac12 W(0)$.

We  now turn to showing that (\textbf{S1}) implies (\textbf{S2}). 
Let us recall the definition of L\'evy--Prokhorov metric, which metrizes the weak convergence of probability measures: Given probability measures $\nu$ and $\sigma$
\[ d_{LP}(\nu, \sigma) = \inf \{ \veps >0 \::\: (\forall A - \te{Borel}) \;\; \nu(A) \leq \sigma(A+ \veps) +\veps 
\te{ and }  \sigma(A) \leq \nu(A+ \veps) +\veps \} \]
where $A+\veps = \{x: d(x,A) < \veps\}$. 

For a given measure $\mu$, we first show that it can be approximated in the L\'evy--Prokhorov metric by an empirical measure of a finite set with arbitrarily many points. That is,  we show that for any $\veps>0$ and any 
$n_0$ there exists $n \geq n_0$ and a set of distinct points $X=\{x_1, \dots, x_n\}$ such that the corresponding empirical measure $\mu_X = \frac{1}{n} \sum_{j=1}^n \delta_{x_j}$ satisfies
$d_{LP}(\mu_X, \mu) < \veps$. 

Let $\veps >0$. We can assume that $\veps < \frac12$.  There exists $R>0$ such that for $Q_R = [-R,R]^N$, $\mu_X(\R^N \backslash Q_R) < \frac{\veps}{2}$. 
For integer $l$ such that  $\sqrt{N} \frac{2R}{l} < \veps$
divide $Q_R$ into $l^N$ disjoint cubes $Q_i$, $i=1, \dots, l^N$ with sides of length $2R/l$. While cubes have the same interiors, they are not required to be identical, namely some may contain different parts of their boundaries, as needed to make them disjoint. Note that the diameter of each cube, $\sqrt{N} \frac{2R}{l}$, is less than $\veps$.
Let $n > n_0$ be such that $ \frac{l^N}{n} < \frac{\veps}{2}$. Let $p = \frac{1}{n}$.
For $i=1, \dots, l^N$ let $p_i = \mu(Q_i)$, $n_i = \lfloor p_i n \rfloor$, and $q_i = n_i p$. 
Note that $0 \leq p_i - q_i \leq p$ and thus $s_q = \sum_{i} q_i \geq \sum_i p_i - l^N p > 1 - \frac{\veps}{2} 
$.  
In each cube $Q_i$ place $n_i$ distinct points and let $\tilde X$ be the set of all such points. 
Note that $\tilde n = \sum_i n_i = s_q n > (1-\veps) n$. Let $\hat X$ be an arbitrary  set of $n - \tilde n$ distinct points in $Q_{2R} \backslash Q_R$. Let $X = \tilde X \cup \hat X$. Note that $X$ is a set of $n$ distinct points. 
Then for any Borel set $A$
\[
\mu(A)  \leq \sum_{i \::\: \mu(A \cap Q_i)>0} \mu(Q_i) + \frac{\veps}{2} 
 \leq \sum_{i \::\: \mu(A \cap Q_i)>0}  \left( \mu_X(Q_i) + p \right) +  \frac{\veps}{2} 
  \leq \mu_X(A+\veps) + \veps.
\]
Similarly
\[ \mu_X(A) \leq \mu(A+ \veps) + \veps. \]
Therefore $d_{LP}(\mu, \mu_X) \leq \veps$.

Consequently there exists a sequence of sets $X_m$ with $n(m)$ points satisfying $n(m) \to \infty$  as $m \to \infty$ for which the empirical measure $\mu_m= \mu_{X_m}$ converges weakly $\mu_m \rhu \mu$ as $m \to \infty$. By assumption \textbf{(S1)}
\begin{equation*}
\iint_{x \neq y} W(x-y) d\mu_{m}(x) d\mu_{X_m} (y) \geq - \frac{1}{n(m)}  B. 
\end{equation*}
Let us first consider the case that $W$ is an upper-semicontinuous function bounded from above.
It follows from Lemma \ref{LSCE} that the energy $E$ is an upper-semicontinuous functional. Therefore
\[ E(\mu) \geq \limsup_{m \to \infty} E(\mu_{m}) \geq \limsup_{m \to \infty} - \frac{1}{n(m)} (B - W(0)) = 0 \]
as desired. 

If $W$ is an upper-semicontinuous function such that there exists $\overline R$ such that $W$ is nondecreasing on $[\overline R, \infty)$ we first note that we can assume that $W(r) \to \infty$ as $r \to \infty$, since otherwise $W$ is bounded from above which is covered by the case above. 
If $\mu$ is a compactly supported probability measure then there exists $L$ such that for all $m$, $\supp \mu_m \subseteq [-L,L]^N$.
Since $W$ is upper-semicontinuous it is bounded from above on compact sets and thus 
upper-semicontinuity of the energy holds. That is $E(\mu) \geq \limsup_{m \to \infty} E(\mu_{m}) \geq 0$ as before.

If $\mu$ is not compactly supported it suffices to show that there exists a compactly supported measure $\tilde \mu$ such that $E(\mu) \geq E(\tilde \mu)$, since by above we know that $E(\tilde \mu) \geq 0$. Note that since $E(\frac12 (\delta_x + \delta_0)) \geq 0$,
$W(|x|) \geq - W(0)$. Therefore  $W$ is bounded from below by $-W(0)$ and $W(0) \geq 0$.

Since $W(r) \to \infty$ as $r \to \infty$ there exists $R_1 \geq \overline R$ such that $W(R_1) \geq \max\{1, \max_{r \leq R_1} W(r) \}$ and $m_1 = \mu(\overline B_{R_1}(0)) > \frac{7}{8}$. Let $R_2$ be such that $W(R_2) > 2 W(R_1)$, and define the constants $m_2 = \mu(\overline B_{R_2}(0) \backslash \overline B_{R_1}(0))$ and $m_3 = \mu(\R^N \backslash \overline B_{R_2}(0))$. Note that $m_1+m_2+m_3=1$. Consider the mapping
\[ P(x) = \begin{cases}
   x & \te{if } |x| \leq R_2 \\
   0 & \te{if } |x|  > R_2.
\end{cases}
\] 
Let $\tilde \mu = P_\sharp \mu$.
Estimating the interaction of particles between the regions provides:
\begin{align*}
 E(\tilde \mu) & \leq E(\mu) + 2 W(0) m_3^2 + 2( W(R_2) + W(0))  m_2 m_3 - 2 (W(R_2) - W(R_1)) m_1 m_3 \\
& \leq E(\mu) + W(R_2) m_3 (m_3 + 4m_2 - m_1) < E(\mu). 
\end{align*}
\end{proof}

\bigskip

As we showed in Theorem \ref{thm:existence-vanish} the property ({\bf HE}) is necessary and sufficient for the existence of ground states when $E$ is defined via an interaction potential satisfying ({\bf H1}), ({\bf H2}) and ({\bf H3b}). The property ({\bf HE}) is posed as a condition directly on the energy $E$, and can be difficult to verify for a given $W$. It is then natural to ask what conditions the interaction potential $W$ needs to satisfy so that the energy $E$ has the property ({\bf HE}). In other words, how can one characterize interaction potentials $w$ for which $E$ admits a global minimizer? 

We do not  address that question in detail, but just comment on the partial results established in the context of $H$-stability of statistical mechanics and how they apply to the minimization of the nonlocal-interaction energy.

Perhaps the first condition which appeared in the statistical mechanics literature states that absolutely integrable potentials which integrate to a negative number over the ambient space are not stable (cf. \cite[Theorem 2]{Dob64} or \cite[Proposition 3.2.4]{Ruelle}). In our language these results translate to the following proposition.

\bprop\label{prop:negative-integral} Consider an interaction potential $w(x)=W(|x|)$ 
where $W$ satisfies the hypotheses ({\bf H1}), ({\bf H2}) and ({\bf H3b}). If $w$
is absolutely integrable on $\rn$ and
	\[
		\int_{\rn} W(|x|)\,dx < 0,
	\]
then the energy $E$ defined by \eqref{eqn:energy} satisfies the condition ({\bf HE}).
\eprop

\begin{proof}
By rescaling, we can assume that  $\int_{\rn} W(|x|)\,dx = -1$. Let $M = \int_{\rn} |W(x)| \, dx$.
Let $R$ be such that $\int_{B_R(0)} W(|x|)\,dx < - \frac{3}{4}$ and $\int_{B_R(0)^c} |W(|x|)| \,dx < \frac{1}{4}$. Consider $n$ large, to be set later, and let
 $\rho(x):=\frac{1}{\omega_N (nR)^N}\chi_{B_{nR}(0)}(x)$, i.e., the scaled characteristic function of the ball of radius $nR$. 
Using the fact that $B_R(0) \subset B_{nR}(y)$ for $|y| < (n-1) R$, we obtain
\begin{align*}
\omega_N^2 (nR)^{2N} \, E(\rho) = &\int_{B_{nR}(0)}\!\int_{B_{nR}(0)} W(|x-y|)\,dxdy \\
						  = &\int_{B_{nR}(0)}\left(\int_{B_{nR}(y)} W(|x|)\,dx\right)\,dy \\
						   \leq & \int_{B_{(n-1)R}(0)} \left( \int_{B_{R}(0)} W(|x|)\,dx + \frac{1}{4} \right)\,dy  +  \int_{B_{nR}(0) \backslash B_{(n-1)R}(0)} M \,dy \\
						   \leq & - \frac12 (n-1)^N R^N \omega_N + N \omega_N n^{N-1}  R^{N} M <0
\end{align*}
if $n$ is large enough.
This shows that the energy $E$ satisfies ({\bf HE}).
\end{proof}

\medskip

An alternative condition for instability of interaction potentials is given in \cite[Section II]{CaSi69}. This condition, which we state and prove in the following proposition, extends the result of Proposition \ref{prop:negative-integral} to interaction potentials which are not necessarily absolutely integrable.

\bprop\label{prop:weighted-int}
Suppose the interaction potential $W$ satisfies the hypotheses ({\bf H1}), ({\bf H2}) and ({\bf H3b}). If there exists $p \geq 0$ for which
		\beqn
			\int_{\R^N} W(|x|)\,e^{-p^2|x|^2}\,dx < 0,
			\label{eqn:weighted-min}
		\eeqn
	then the energy $E$ defined by \eqref{eqn:energy} satisfies the condition ({\bf HE}).
\eprop

\begin{proof}
	Let $p\geq 0$ be given such that the inequality \eqref{eqn:weighted-min} holds. Since the case $p=0$ has been considered in Proposition \ref{prop:negative-integral}, we can assume $p>0$.	
Consider the function
		\[
			\rho(x) = \frac{p^N}{\pi^{N/2}}\,e^{-2p^2|x|^2}.
		\]
Clearly $\rho\in L^1(\rn)$ and $\|\rho\|_{L^1(\rn)}=1$; hence, it defines a probability measure on $\rn$. Consider the linear transformation on $\mathbb{R}^{2N}$ given by
	\[
		u=x-y,\qquad v=x+y.
	\]
We note that the Jacobian of the transformation is 2. Thus
	\beqns
		\begin{aligned}
			E(\rho) &= \int_{\rn}\!\int_{\rn} W(|x-y|)\,e^{-2p^2|x|^2}\,e^{-2p^2|y|^2}\,dxdy \\
							&= \frac12 \int_{\rn}\!\int_{\rn} W(|u|)\,e^{-p^2|u+v|^2/2}\,e^{-p^2|u-v|^2/2}\,dudv \\
							&= \frac12 \int_{\rn}\!\int_{\rn} W(|u|)\,e^{-p^2(|u|^2+|v|^2)}\,dudv \\
							&= \frac12 \int_{\rn}\left(\int_{\rn}W(|u|)\,e^{-p^2|u|^2}\,du\right)e^{-p^2|v|^2}\,dv  < 0.
		\end{aligned}
	\eeqns
Hence, the energy $E$ satisfies ({\bf HE}).
\end{proof}

\begin{remark}
Another useful criterion can be obtained by using the Fourier transform, as also noted in \cite{Ruelle}. Namely if $w \in L^2(\R^N)$,
for measure $\mu$ that has a density $\rho \in L^2(\R^N)$, 
by Plancharel's theorem
\[ E(\mu) =  \int_{\rn}\!\int_{\rn} w(x-y)\,d\mu(x)d\mu(y) = \int_{\R^N} \hat w(\xi) |\hat \rho(\xi)|^2 d \xi. \]
So if the real part of $\hat w$ is positive, the energy does not have a minimizer.

This criterion can be refined. By Bochner's theorem the Fourier transforms of finite nonnegative measures are precisely the positive definite functions. Thus we know which family of functions, $\hat \rho$ belongs to. Hence we can formulate the following criterion:
If $w \in L^2(\R^N)$ and there exists a positive definite complex valued function $\psi$  such that $ \int \hat w(\xi) | \psi^2(\xi)| d \xi \leq 0$ then the energy $E$ satisfies the condition  ({\bf HE}).
\end{remark}

\bigskip
\noindent {\bf Acknowledgments.}
 The authors would like to thank the Center for Nonlinear Analysis of the Carnegie Mellon University for its support, and hospitality during IT's visit. 
RS was supported by  the Funda\c{c}\~{a}o para a Ci\^{e}ncia e a Tecnologia (Portuguese Foundation for Science and Technology) through the
Carnegie Mellon Portugal Program
under Grant SFRH/BD/33778/2009.
 DS is grateful to  NSF (grant DMS-1211760) and FCT (grant UTA CMU/MAT/0007/2009).
 IT was also partially supported by the Applied Mathematics Laboratory of the Centre de Recherches Math\'{e}matiques.  The research was also supported by NSF PIRE grant  OISE-0967140.

\bibliography{lit}{}
\bibliographystyle{plain}	

\end{document}